\documentclass[12pt]{article}
\usepackage{amstex,amscd}
\makeatletter
\def\section{\@startsection {section}{1}{\z@}{3.5ex plus1ex minus
    .2ex}{2.3ex plus.2ex}{\Large\bf}}
\def\subsection{\@startsection{subsection}{2}{\z@}{3.25ex plus1ex
    minus.2ex}{1.5ex plus.2ex}{\reset@font\large\bf}}
\def\subsubsection{\@startsection{subsubsection}{3}{\z@}{3.25ex plus
    1ex minus.2ex}{1.5ex plus.2ex}{\reset@font\normalsize\bf}}
\def\@aftersepkern{\kern-.5em}
\def\postchapter{.}
\def\postsection{.\@aftersepkern}
\def\postsubsection{.\@aftersepkern}
\def\postsubsubsection{.\@aftersepkern}
\def\postparagraph{.\@aftersepkern}
\def\postsubparagraph{.\@aftersepkern}
\newbox\znak
\setbox\znak=\hbox{{\Large\bf \S\,}}

\def\@sect#1#2#3#4#5#6[#7]#8{\ifnum #2>\c@secnumdepth
     \let\@svsec\@empty\else
     \refstepcounter{#1}\edef\@svsec{%
     \csname pre#1\endcsname
     \csname the#1\endcsname
     \csname post#1\endcsname
\hskip 1em}\fi
     \@tempskipa #5\relax
      \ifdim \@tempskipa>\z@
        \begingroup #6\relax
          \@hangfrom{\hskip #3\relax\@svsec
                   }{\interlinepenalty \@M \ignorespaces#8\par}%
        \endgroup
       \csname #1mark\endcsname{\protect\ignorespaces #7}\addcontentsline
         {toc}{#1}{\ifnum #2>\c@secnumdepth \else
                      \protect\numberline{\csname the#1\endcsname
                     \csname post#1\endcsname
}\fi
\protect\ignorespaces #7}\else
        \def\@svsechd{#6\hskip #3\relax  
                   \@svsec #8\csname #1mark\endcsname
                      {\protect\ignorespaces #7}\addcontentsline
                           {toc}{#1}{\ifnum #2>\c@secnumdepth \else
                             \protect\numberline{\csname the#1\endcsname
                          \csname post#1\endcsname
}\fi
\protect\ignorespaces #7}}\fi
     \@xsect{#5}}
\def\@ssect#1#2#3#4#5{\@tempskipa #3\relax
   \ifdim \@tempskipa>\z@
     \begingroup #4\@hangfrom{\hskip #1}{\interlinepenalty \@M
\ignorespaces #5\par}\endgroup
   \else \def\@svsechd{#4\hskip #1\relax
\ignorespaces #5}\fi
    \@xsect{#3}}
\def\@makechapterhead#1{%
  \vspace*{50\p@}%
  {\parindent \z@\raggedright
    \ifnum \c@secnumdepth >\m@ne
     \huge\bf \@chapapp{} \thechapter\postchapter
    \par
    \vskip 20\p@ \fi
    \Huge \bf
    \ignorespaces #1\par
    \nobreak
    \vskip 40\p@
  }}
\def\@makeschapterhead#1{%
  \vspace*{50\p@}%
  {\parindent \z@ \raggedright
    \Huge \bf
    \ignorespaces #1\par
    \nobreak
    \vskip 40\p@
  }}
\def\@begintheorem#1#2{\trivlist \item[\hskip \labelsep{\sc #1\ #2.}]\sl}
\def\@opargbegintheorem#1#2#3{\trivlist
      \item[\hskip \labelsep{\sc #1\ #2\ (#3).}]\sl}
\def\@endtheorem{\endtrivlist}
\newtheorem{theorem}{Theorem}[section]
\newtheorem{proposition}{Proposition}[section]
\newtheorem{lemma}{Lemma}[section]
\newcounter{definition}[section]

\newenvironment{definition}{\par\medskip\noindent\refstepcounter{definition}%
           {\sc Definition \arabic{section}.\arabic{definition}.} }%
           {\medskip\par}
\newcounter{example}[section]

\newenvironment{example}{\par\medskip\noindent\refstepcounter{example}%
           {\sc Example \arabic{section}.\arabic{example}.} }%
           {\medskip\par}
\newenvironment{indentedtext}{\par\noindent\advance\leftskip by1cm}{\par}
\makeatother
\author{Theodore Voronov}

\title{Supermanifold Forms and Integration. A Dual Theory.
\thanks{The research was partially supported by the
Russian Foundation for Basic Research
(grant No 94-01-01444a)}}

\date{{\small Department of Mathematics, Moscow State University\\
       Vorobjovy Gory, Moscow 119899, Russia\\
    {\it E-mail:\/} {theodore$@$mech.math.msu.su, theodore$@$nw.math.msu.su}}}

\numberwithin{equation}{section}
\DeclareMathSymbol\leqslant{\mathrel}{AMSa}{"36} 
\DeclareMathSymbol\geqslant{\mathrel}{AMSa}{"3E} 
\renewcommand{\ge}{\geqslant}
\renewcommand{\le}{\leqslant}

\begin{document}
\maketitle
%
%
\renewcommand{\a}{\alpha    }
\renewcommand{\b}{\beta    }
\renewcommand{\c   }{\gamma    }
\newcommand{\C   }{\Gamma    }
\renewcommand{\d   }{\delta    }
\newcommand{\D   }{\Delta    }
\newcommand{\e   }{\varepsilon}
\newcommand{\g   }{\gamma    }
\newcommand{\G   }{\Gamma    }
\newcommand{\h   }{\eta    }
\renewcommand{\l   }{\lambda    }
\renewcommand{\L   }{\Lambda    }
\newcommand{\m   }{\mu    }
\newcommand{\n   }{\nu    }
\renewcommand{\O   }{\Omega }
\renewcommand{\o   }{\omega    }
\newcommand{\p   }{\pi    }
\renewcommand{\P   }{\Pi }
\newcommand{\s   }{\sigma    }
\newcommand{\Ss   }{\Sigma    }
\renewcommand{\t   }{\tau    }
\renewcommand{\u   }{\upsilon    }
\newcommand{\U   }{\Upsilon    }
\renewcommand{\v   }{\upsilon    }
\newcommand{\V   }{\Upsilon    }
\newcommand{\x   }{\xi    }
\newcommand{\X   }{\Xi    }
\newcommand{\y   }{\upsilon    }
\newcommand{\Y   }{\Upsilon    }
\newcommand{\z   }{\zeta    }
\newcommand{\Z   }{\Zeta    }
\newcommand{\ps   }{\psi    }
\newcommand{\ph   }{\varphi    }
\newcommand{\PH   }{\varPhi    }
%
%
%
%
%
%
\newcommand{\mqquad}{{\!\!\!\!\!\!\!\!\!\!\!\!}} 
\newcommand{\mquad}{{\!\!\!\!\!\!}}              
\newcommand{\qed}{{$\quad\Box$}}

\newcommand{\gl}{{\operatorname{GL}}}
\newcommand{\mat}{{\operatorname{Mat}}}
\newcommand{\glpq}{{\gl(p|q)}}
\newcommand{\ber}{{\operatorname{Ber}\,}}
\renewcommand{\div}{{\operatorname{div}\,}}
\newcommand{\sign}{{\operatorname{sign}}}
\newcommand{\rank}{{\operatorname{rank}\,}}
\newcommand{\codeg}{{\operatorname{codeg}\,}}
\newcommand{\addeg}{{\operatorname{add.deg}\,}}

\newcommand{\der}[2]{{\frac{\partial #1}{\partial #2}}}
\newcommand{\derl}[2]{{{\partial #1}/{\partial #2}}}
\newcommand{\secder}[3]{{\frac{\partial^2 #1}{{\partial #2}{\partial #3}}}}
\newcommand{\secderl}[3]{{{\partial^2 #1}/{\partial #2}{\partial #3}}}
\newcommand{\thirder}[4]{{\frac{\partial^3 #1}{{\partial #2}{\partial #3}%
                               {\partial #4}}}}

\renewcommand{\L}{{\cal L}}
\renewcommand{\D}{{\,\,}{\bar{\smash{\!\!d}}}}
\newcommand{\DD}{{}{\bar{\smash{\d}}}}

\newcommand{\pak}{{p_A{}^K}}
\newcommand{\pal}{{p_A{}^L}}
\newcommand{\pag}{{p_A{}^G}}
\newcommand{\pbk}{{p_B{}^K}}
\newcommand{\pbl}{{p_B{}^L}}
\newcommand{\wfk}{{w_F{}^K}}
\newcommand{\wfl}{{w_F{}^L}}
\newcommand{\wfg}{{w_F{}^G}}
\newcommand{\wgk}{{w_G{}^K}}
\newcommand{\wgl}{{w_G{}^L}}
\newcommand{\wk}{{w^K}}
\newcommand{\va}{{v_A}}
\newcommand{\fk}{{f^K}}
\newcommand{\yk}{{Y^K}}
\newcommand{\zlk}{{Z_L{}^K }}
\newcommand{\dkl}{{\d_K{}^L}}
\newcommand{\dlk}{{\d_L{}^K}}
\newcommand{\dfg}{{\d_F{}^G}}
\newcommand{\ddfk}{{\der{\d}{f^K}}}

\newcommand{\kt}{{\tilde K}}
\newcommand{\lt}{{\tilde L}}
\newcommand{\at}{{\tilde A}}
\newcommand{\bt}{{\tilde B}}
\newcommand{\dt}{{\tilde \d}}
\newcommand{\ft}{{\tilde F}}
\newcommand{\gt}{{\tilde G}}

\newcommand{\mbig}{M^{n|m}\times{\Bbb R}^{r|s}}
\newcommand{\ubig}{U^{n|m}\times{\Bbb R}^{r|s}}

\newcommand{\pone}{{p_1}}
\newcommand{\ptwo}{{p_2}}
\newcommand{\woo}{{w_{11}}}
\newcommand{\wot}{{w_{12}}}
\newcommand{\wto}{{w_{21}}}
\newcommand{\wtt}{{w_{22}}}

\newcommand{\om}[4]{{\O^{#1|#2}_{#3|#4}}}
\newcommand{\Om}[2]{{\boldsymbol\O^{#1|#2}}}
\newcommand{\si}[4]{{\Sigma^{#1|#2}_{#3|#4}}}
\newcommand{\Si}[2]{{\boldsymbol\Sigma^{#1|#2}}}

\newcommand{\dotarrow}[8]{\begin{picture}(0,0) 
                            \multiput(0,0)(#4,#5){#6}
                               {\line(#1,#2){#3}}    
                            \put(#7,#8){\vector(#1,#2){#3}} 
                         \end{picture}}
\newcommand{\dar}{\dotarrow{4}{3}{10}{12}{9}{3}{36}{27}}
\newcommand{\larrow}[2]{\begin{picture}(0,0)
                         \put(#1,#2){\begin{picture}(0,0)
                                       \put(0,0){\dar}
                                       \put(15,19){${\scriptstyle \h}$}
                                     \end{picture}}
                        \end{picture}}
\newcommand{\rarrow}[2]{\begin{picture}(0,0)
                         \put(#1,#2){\begin{picture}(0,0)
                                       \put(0,0){\dar}
                                       \put(15,19){${\scriptstyle \iota}$}
                                     \end{picture}}
                        \end{picture}}


\begin{abstract}
      We investigate forms on supermanifolds defined as Lagrangians of
``copaths''. For this, we consider direct products
$M^{n|m}\times\Bbb R^{r|s}$
and study isomorphisms corresponding to simultaneously advancing
the number of additional parameters  $r|s$ and the number of  equations.
We define an exteriour differential in terms of variational derivatives
w.r.t. a copath and establish its main properties. In the resulting
stable picture we obtain infinite complexes
$\D:\Om{r}{s}\to\Om{r+1}{s}$ for $M^{n|m}$,
where $0\le s\le m$ and $r$ can be any integer. For
$r\ge 0$ a canonical isomorphism with forms constructed as
Lagrangians of $r|s$-paths is established.
We discover the ``lacking half'' of forms on supermanifolds:
$r|s$-forms with
$r<0$, previously unknown except for $s=m$.
(They have been partly replaced earlier by an augmentation of
the ``non-negative'' part of the complexes.) 
The study of these questions is in progress now.
\end{abstract}

%
%
%
%
%
\section*{Introduction}
Most of ``supermathematics'' can be obtained more or less easily by
extending ``purely even'' constructions with the help of the sign rule.
Supermanifold integration is not the case. Since
the times when
the Berezin integral was discovered and the notion of supermanifold came
into being, an adequate theory of ``super'' forms suiting for integration
has been a great puzzle. (The naive definition of differential forms
in super case has nothing to do with the Berezin integration.)
This was not for chance (see~\cite{git} for discussion).
The construction
of the theory in question took great effort in years 1975--1985.
See Bernstein and
Leites~\cite{ber1,ber2}, Gajduk, Khudaverdian and
Schwarz~\cite{gajduk},
Baranov and Schwarz~\cite{barsh}, Voronov and
Zorich~\cite{first,our,bord,dan},
Voronov~\cite{dis,git}, and others. An essencial feature of the
theory of forms developped in~\cite{first,our,bord,%
dan,dis,git} is that forms of (even) degree greater than the even
dimension of a supermanifold are present necessarily.

In this paper we develop a ``dual'' approach, 
in which forms are
treated as Lagrangians of ``copaths'' (systems of equations, which
may specify nonsingular submanifolds or may not). To introduce a
differential, we extend supermanifolds in consideration
by introducing additional variables and
study a sort of ``stabilization'', reminiscent of  $K$-theoretic methods.
This leads us to a discovery of a ``lacking half'' of the previously
developped
theory of forms, namely,
forms (and perhaps cohomology) of negative degree.
All results are new.
The study of these
questions is in progress now.

(I would like to note that the ``dual approach'' to supermanifold
integration has been
independently investigated by O.M.Khudaverdian who obtained
important results.)

The paper is organized as follows. In \S 1 we present a survey of forms as
Lagrangians of paths (a ``standard'' theory).
For greater details we refer to~\cite{git}. In \S 2 we
introduce Lagrangians of copaths and study the variation of the
corresponding action. In \S 3 we study ``mixed Lagrangians'' and prove the
``stabilization theorem'' (Theorem~\ref{iso}). We
introduce a differential for ``mixed forms''
and establish its main properties
(Theorem~\ref{main}). It is shown how
the complexes of ``standard'' forms and dual forms are sewed
together in the ``stable'' picture.

Throughout the paper we use the notation of the book~\cite{git}.
I wish to thank J.N.Bernstein, O.M.Khudaverdian and M.A.Shubin. Dis\-cus\-sions
of super\-manifold forms and integration with them has stimulated this
research greatly and I am much grateful to them.
%
%
%
%
%
%
%
\section{Forms as Lagrangians of multidimensional paths.}

\newcommand{\irs}{{I^{r|s}}}
\newcommand{\xdot}{\dot x}

Let $M=M^{n|m}$ be a $n|m$-dimensional supermanifold. An {\it $r|s$-path\/}
in $M$ is a smooth map $\G:\irs\to M$, where $\irs=I^r\times\Bbb
R^{0|s}\subset\Bbb R^{r|s}$ is taken with a fixed boundary. (Fixing
boundary is necessary for Berezin integration over a domain, see~\cite{git}
and also below.)
Denote coordinates on $M$ by $x^A$, on $\irs$ by $t^F$, so the path $\G$
is presented as $x^A=x^A(t)$.

Recall that a Lagrangian function or simply a Lagrangian on $M$ is a
function on the tangent bundle $TM$. It depends on a position and on a
single tangent vector. Integrating a Lagrangian over a path (that is, a
$1|0$-path), one obtains the action of the path. Note that in this case
$dx/dt$ is an even vector.

We shall call a function of a position and of an array of $r+s$
tangent vectors, of which $r$ are even and $s$ are odd,
a {\it (generalized) Lagrangian on \/} $M$. We omit the word
``generalized'' in the following. A Lagrangian $L=L(x,\dot x)$,
where $\dot x=(\dot
x_F{}^A)$, defines an action for $r|s$-path:
\begin{equation}\label{sg}
    S[\G]=\int\limits_{\irs}\!Dt\,L\left(x(t),\der{x}{t}(t)\right).
\end{equation}
Here $Dt$ stands for a Berezin ``volume element'' on superspace.

Recall some facts about the Berezin integral (see~\cite{git}).
Roughly
speaking, Berezin integration consists in singling out the top-order
coefficient in the power expansion in  odd variables and in usual
integration over even variables.
For bounded domains
a step-function $\theta(g)$ is inserted. Here the boundary is specified
by the equation $g=0$. Since the even function $g$ may contain nilpotents,
integration uses more data than just the underlying
set of points of the domain.
The formula for the change of variables in the usual Riemann or Lebesgue
integral contains the absolute value of the Jacobian. Its analogue for the
change of variables in the Berezin integral is
the function $\ber_{1,0}\,J=\ber\,J\cdot\sign\det J_{00}$.
Here $J$ is the Jacobian matrix of the transformation of variables and
$J_{00}$ stands for its even-even block. For a block square matrix
(blocks are specified by attributing ``parity'' to columns and rows)
with even elements in the diagonal blocks and odd elements in the
antidiagonal blocks
the Berezinian is defined as
\begin{equation}
  \begin{aligned}
    \ber J&=\det(J_{00}-J_{01}{J_{11}}^{-1}J_{10})\cdot(\det J_{11})^{-1}\\
          &=\det J_{00}\cdot(\det(J_{11}-J_{10}{J_{00}}^{-1}J_{01}))^{-1}.
  \end{aligned}
\end{equation}
It is  a unique (in essence, up to trivial changes) function on
invertible matrices near unity that is multiplicative.
It is essential that the
general linear group in supercase possesses four, not two, connected
components. They are specified by various combinations of the signs of the
even-even and the odd-odd blocks.
Two functions, $\det J$ and $|\det J|$, in supercase
turn into four: $\ber_{\!\a\b}\,J=\ber J\cdot(\sign\det J_{00})^{\a}
\cdot(\sign\det J_{11})^{\b}$, where $\a,\b=0,1$.
Recall that an orientation of vector space is a class of frames
modulo linear transformations with positive determinant. We see that in
supercase principally different notions of orientation are possible. The
most important are  a {\it $(+|-)$-orientation\/} preserved by
linear transformations with $\det J_{00}>0$ and
a {\it $(-|-)$-orientation\/}
preserved by transformations with $\ber J_{00}>0$ (see~\cite{git}). 
We also arrive to four possible types of ``volume elements'' $D_{\!\a\b}$:
\begin{equation}
    D_{\!\a\b}x=D_{\!\a\b}x'\cdot\frac{D_{\!\a\b}\,x}{D_{\!\a\b}\,x'},
    \text{\quad where\quad}\frac{D_{\!\a\b}x}{D_{\!\a\b}x'}=
    \ber_{\!\a\b}\,\left(\der{x}{x'}\right).
\end{equation}
Here $Dx=D_{0,0}x$ and $D_{0,1}x$ are the analogues of what is called the
``oriented volume element'' in the classical case and $D_{1,0}x$ and
$D_{1,1}x$ are the analogues of the ``unoriented volume element'' (denoted
sometimes by $|dx|$ or similar). Actually, an {\it unoriented volume
element\/} in supercase is $D_{1,0}x$. The expression $f(x)D_{1,0}x$,
where $f$ is a function, can be integrated over supermanifold without any
conditions on orientation. At the same time, integration of the volume form
like $f(x)Dx$ requires that the supermanifold in question should be endowed
with a  $(+|-)$-orientation.

Return to the integral~(\ref{sg}). Consider changes of parametrization for
the path $\G$.
\begin{proposition}
    The action $S[\G]$ is invariant under arbitrary reparametrization iff
\begin{equation}
    L(x,h\xdot)=\ber_{\!1,0}\, h\cdot L(x,\xdot)
\end{equation}
for any $h\in\gl(r|s)$.
\end{proposition}
The Lagrangians with such property are the multidimensional analogue of the
homogeneous Lagrangian functions (in the classical sense) like
the Lagrangian
of the free particle taken in the form $m|\xdot|$.
They are well-defined only if all vectors $\xdot_F$ are linear independent,
in particular no vector can vanish. That means that dealing with such
Lagrangians one should consider only paths which are immersions.
This would be a bad starting point for a theory of forms.
Instead  we need to restrict the
possible reparametrization of paths by certain orientation conditions.
\begin{definition}\label{cov1}
    A Lagrangian L is called {\it covariant (of the first kind)\/}, if
\begin{equation}
    L(x,h\xdot)=\ber h\cdot L(x,\xdot),
\end{equation}
and it is called {\it covariant (of the second kind)\/}, if
\begin{equation}
    L(x,h\xdot)=\ber_{\!0,1}\, h\cdot L(x,\xdot),
\end{equation}
for any $h\in \gl(r|s)$.
\end{definition}

It follows that for covariant Lagrangians of the first kind the action
is invariant under any changes of parameters that preserve
the $(+|-)$-orientation of the path $\G$. The same is true for covariant
Lagrangians of the second kind and the $(-|-)$-orientation.
The natural widest domain of definition for covariant Lagrangians
considered on some $U\subset M$ is
an open domain $W_{r|s}(U)\subset \underbrace{TU\oplus\ldots\oplus TU}_{r}
\oplus\underbrace{\P TU\oplus\ldots\oplus\P TU}_{s}$ specified by the
condition \quad $\rank (\xdot_{\a}^{\m})=s$.
(Here $\a=1\ldots s$ runs over odd
indices $F$ and $\m=1\ldots m$ runs over odd $A$.) By $\P$ we denote the
parity reversion functor. We consider
covariant Lagrangians defined everywhere on $W_{r|s}$.

\begin{definition}
    A smooth map of supermanifolds is {\it proper\/}, if it induces
monomorphisms of the odd subspaces of the tangent spaces (at each point).
\end{definition}

Obviously, the pull-back of covariant Lagrangians is well-defined provided
the map is proper.

There exist nonvanishing covariant Lagrangians with $r>n$, if $m>0$.
This will be obvious if one recalls the definition of the ``naive''
differential forms in supercase. The algebra of such forms is generated by
differentials $dx^A$, which have parity opposite to that of coordinates.
So if $m>0$, there exist products of $d\x^{\m}$ of arbitrary high degree.
Since any ``naive'' differential form can be integrated over $r|0$-paths,
this provides the desired examples of covariant Lagrangians.
Examples with $s>0$ exist too (see~\cite{git}).

Consider the variational derivative of the action ~(\ref{sg}).
\begin{proposition}
\begin{multline}\label{var1}
      \frac{\d S}{\d x^A}=\der{L}{x^A}-(-1)^{\at\ft}\der{}{t^F}
                            \left(\der{L}{\xdot_F{}^A}\right) \\
       =\der{L}{x^A}-(-1)^{\at\ft}\xdot_F{}^B\secder{L}{x^B}{\xdot_F{}^A}\\
          -(-1)^{\at\ft}\frac{1}{2}\ddot x_{FG}^B
\left(\secder{L}{\xdot_G{}^B}{\xdot_F{}^A}+(-1)^{\gt\ft+\at(\gt+\ft)}
\secder{L}{\xdot_F{}^B}{\xdot_G{}^A}\right).
\end{multline}
\end{proposition}
\begin{definition}\label{form}
    A covariant Lagrangian $L$ is called a {\it form\/} (of {\it degree\/}
$r|s$), if the last term in~(\ref{var1}) identically vanish:
\begin{equation}\label{bas1}
\secder{L}{\xdot_F{}^A}{\xdot_G{}^B}+(-1)^{\ft\gt+\bt(\ft+\gt)}
\secder{L}{\xdot_G{}^A}{\xdot_F{}^B}=0
\end{equation}
for any $A,B,F,G$.
\end{definition}
The equations ~(\ref{bas1}) guarantee that the variation does not depend on
the ``acceleration''. To proceed with the physical analogy, such Lagrangians
correspond to interaction terms like $-e A_a\xdot^a$, for a charged particle
in the electromagnetic field. One can check that in purely even case the
equations~(\ref{bas1}) together with the covariance condition are
equivalent to that $L$ is multilinear and skew-symmetric function of
tangent vectors.
(In general, the equations~(\ref{bas1}) imply that $L$
is skew-symmetric and multiaffine both in even rows and even columns of
$(\xdot_F{}^A)$.)
So in this case our definition reproduces the classical
concept of an exteriour form. The equations~(\ref{bas1}) are functorial,
so forms can be pulled back along proper morphisms.

Forms can by no means be multilinear in odd vectors, since they should be
homogeneous of degree $-1$ w.r.t. them. The equations~(\ref{bas1}) could be
regarded as nontrivial substitution for multilinearity.
The odd-odd part of ~(\ref{bas1}) is a system of partial
differential equations for a function of an even matrix.
A remarkable
fact is that this system (or similar) appears in classical integral geometry,
and is typical in
connection with Radon-like integral transforms
(see~\cite{john,ggg}). There are many unexpected
links of integration on
supermanifolds with classical integral geometry (see
{}~\cite{git,mor,dan,shan}).

\begin{definition}
    A {\it differential} of $r|s$-form $L$ is
an $(r+1\,|\,s)$-form $\D L$ defined
as follows:
\begin{equation}
    \D L=(-1)^r\xdot_{r+1}^A\left(
               \der{L}{x^A}-(-1)^{\at\ft}\xdot_F{}^B
               \secder{L}{x^B}{\xdot_F{}^A} \right).
\end{equation}
\end{definition}

One can check that the operation $\D$ is well-defined. The identity
$\D\,^2=0$ holds. Forms can be integrated over {\it singular manifolds\/}
$f: P^{r|s}\to M^{n|m}$ ($f$ is required to be proper), provided the
orientation of the appropriate kind is fixed on $P^{r|s}$. If $P^{r|s}$
is considered together with a chosen boundary,
then the Stokes formula holds:
\begin{equation}
  \int\limits_{(P^{r|s},\,f)}\!\!\D L\,=\!\!
  \int\limits_{(Q^{r-1|s},\,f\circ i)}\!\!\!\!L,
\end{equation}
where $L$ is a $(r-1\,|\,s)$-form and $i:Q^{r-1|s}\to P^{r|s}$ the
inclusion of the boundary.

The problem of the theory of supermanifold forms, in the sense of
Definition~\ref{form}, is that no explicit description is at hand.
An exception is some special cases. For $s=0$ there is $1-1$-correspondence
with the ``naive'' differential forms. For $s=m$ (odd codimension zero)
it can be proved that $r|m$-forms are in bijection with so called integral
forms of Bernstein and Leites~\cite{ber1}, if $r\ge 0$ (see \cite{git}).
Note that these cases ($s=0$ or $s=m$) are exceptional also because here
the difference between forms of first and second kind is immaterial.

For forms of the second kind there is a construction supplying plenty
of examples. Let $\o=\o(x,dx)$ be a Bernstein--Leites pseudodifferential
form (p.d.f.)~\cite{ber2}.
It is a function that need not be polynomial in $dx^A$.
(For example, $e^{-(d\x)^2}$ is possible.)
This is a beautiful construction but one has to pay for this beauty.
First, there is
no grading. Second, if $\o(x,dx)$ sufficiently rapidly decreases at the
infinity of $d\x^{\m}$, it can be integrated over a singular manifold of
arbitrary dimension ~\cite{ber2}, but {\sl only\/} in the
$(-|-)$-oriented case. It doesn't work for $(+|-)$-orientation.
An integral of $\o$ over $(P,f)$ is defined
as the integral of $f^*\o$.
For a $(-|-)$-oriented supermanifold  $M$ an integral of p.d.f. $\o$
is defined ~\cite{ber2} as:
\begin{equation}
  \int\limits_M\o\: :=\int\limits_{\P TM}\!\! Dx\,D(dx)\,\,\o(x,dx).
\end{equation}
(In the classical case this coincides with extracting the term of
maximal degree from an inhomogeneous form and integrating it over a
manifold.) This leads to the following integral transforms, which
map
Bernstein--Leites p.d.f.'s to $r|s$-forms of the second kind
(see~\cite{barsh,our,git}):
\begin{equation}
  L(x,\xdot)=\int\limits_{\Bbb R^{s|r}}  \! D(dt)\,\,\o(x,dt\cdot\xdot).
\end{equation}
No transform like this is available for forms of the first kind.

One of the motivations for the theory of forms on supermanifolds is the
problem of supermanifold cohomology. The naive differential forms as well
as integral forms and pseudodifferential forms reproduce the cohomology
of the underlying space. There is some hope that forms
in the sense of Definition~\ref{form} can give richer cohomology.
The actual situation is the following. From homotopy point of view, the
category of supermanifolds is equivalent to the category of vector bundles
{}~\cite{git}
(with proper maps and fiberwise monomorphisms, respectively).
So there is a problem
of developping a ``cohomology theory'' for vector bundles (in the indicated
category). There exists a spectral sequence analoguous to the
Atiyah--Hirzebruch sequence~\cite{dan,git}. Its limit term is
adjoined to the cohomology of (global) forms on $M^{n|m}$ and the second
term is the cohomology of the underlying space with coefficients in the
cohomology of forms on $\Bbb R^{0|m}$. Here $\Bbb R^{0|m}$ plays the r\^ole
of ``point'' in standard homotopy theory. A remarkable fact is that
this ``point cohomology'' for some $s\neq 0,m$ is not trivial.
There exist
estimates from the below for it (for details see~\cite{dan,git}).
Unfortunately, the complete calculation for ``supermanifold de~Rham
cohomology'' is not available for the present moment (see~\cite{git} for
discussion).

The results on cohomology described above are seriously based on the
``Cartan calculus'' for $r|s$-forms. It was found ~\cite{git} that the
homotopy identity (or the ``Main Formula of the Differential Calculus for
Forms'')
is not valid unless we introduce in a formal way the ``exact $0|s$-forms''
$B^{0|s}$. Actually they are defined as discrepancies for the generalized
homotopy identity ~\cite{git}. Only after the complexes of $\cdot|s$-forms
are augmented in this way, the cohomology theory can be developped.
At the other hand, the comparison with Bernstein--Leites integral forms,
which are isomorphic to $r|m$-forms, if $r\ge 0$, but which are defined
and nonvanishing for $r<0$, suggests~\cite{git} a search for forms of
negative (even)
degree for an arbitrary $s$.
%
%
%
%
%
%
%
\section{Copaths and dual Lagrangians.}

Consider an open domain $U$ in a supermanifold $M^{n|m}$. A {\it copath} in
$U$ is an array of functions $f^K\in C^\infty{(U)}$ enumerated
by indices $K$,
which run over ``even'' and ``odd'' values so that $f^K$ can be formally
treated as coordinates on some ${\Bbb R}^{p|q}$. The dimension ${p|q}$ is
called a
(formal) {\it codimension} of the copath. A copath can be also treated as a
``framed ideal'' in $C^\infty{(U)}$, i.e. an ideal with  chosen generators.
The corresponding closed subspace of $U$ is a submanifold (of codimension
$p|q$), if the functions $\fk$ are independent.

Consider in a formal way the following integral (the ``action''):
\begin{equation}
   S[f] = \int\limits_{U^{n|m}}\! Dx\, \d(f)\, {\cal L}\left(x,{{\partial f}
                                                \over {\partial x}}\right).
\label{S}
\end{equation}
We leave aside here the questions of the convergence of the
integral, the influence of the boundary and
the existence of the delta-function $\d(f)$.

Here $\L=\L(x,p)$  is a function
of the local coordinates $x^A$ and the components of the covectors
(``momenta'')
$p^K=(\pak)$. (Warning: one should not mix the notation for the momenta
$p=(\pak)$
with $p$ as the number of even functions $f^K$ !)
As a tensor object $\L$ should be a component of an unoriented
volume form, so to make the integral ~(\ref{S})
independent on the choice of coordinates.
We shall call such functions $\L$
{\it (dual) Lagrangians\,} on $M$.

\begin{definition}\label{deg}
The number of $p^K$ and its complement are called
the {\it codegree\,} and the {\it degree\,} of $\L$, respectively:
\begin{align}
    &\codeg\L:=p|q,\\
    &\deg\L:= \dim M-\codeg\L=n-p\,|\,m-q.
\end{align}
\end{definition}

\vspace{-4ex  minus 1ex}
Consider the variation of the action~(\ref{S}).
This needs certain caution. First, one
can treat ~(\ref{S}) as a usual action for the field $\fk$ with a
singular
Lagrangian density $\Lambda(x,f,\,\derl{f}{x})=
\d(f)\,\L(x,\,\derl{f}{x})$.
The formal calculation of the variational derivative
according to the
Euler--Lagrange  formulae leads to
\begin{equation}
\d S=\pm  \int Dx\, \yk \,\left(\ddfk \L -(-1)^{\at\kt+(\at+\kt)\dt}
    \der{}{x^A}\left(\d \der{\L}{\pak}\right)\right), \label{eulagr}
\end{equation}
where $Y=(Y^K)$ is the variation of $f$: $f^K\mapsto f^K+\e\,Y^K$, $\e^2=0$.
(The common sign in~(\ref{eulagr})
depends on dimensions and the parity of $Y$ and is
irrelevant.)
Here the substitution $f^K=f^K(x)$, $\pak=\derl{f^K}{x^A}$
before taking the ``total'' derivative by $x^A$
is supposed, as usual.

One should note, however, that this variation of $S$ corresponds only to a
limited class of variations of $f^K$, namely to the variations compactly
supported in $U$.
If we are going to consider our ``copaths'' as corresponding to surfaces in
$M$, then we should also consider variations which correspond
to the change of
equations $f^K=0$ for an equivalent system. These variations of $f^K$
are of the
form $Y^K=f^L\,Z_L{}^K$, where arbitrary $Z_L{}^K\in C^\infty(U)$
are not supposed
to be compactly
supported.

\begin{proposition}
    The variation of the action is
    \begin{equation}
      \d S=\pm\int Dx \,\zlk \d(f) \left(-(-1)^{\kt}\dkl \L +
                        (-1)^{\at(\kt+\lt)}\pal\der{\L}{\pak}\right)
    \label{frame}
    \end{equation}
in this case.
\end{proposition}
\vspace{-2ex minus 1ex}
{\sc Proof. }Straightforward calculation, using the easily obtained identity
$f^L\derl{\d}{f^K}=
-(-1)^{\kt\lt}\dlk\,\d(f)$.\qed

\medskip
Let the variation ~(\ref{frame}) vanish identically. This implies the
equation
\begin{equation}
 (-1)^{\at(\kt+\lt)}\pal\der{\L}{\pak}=
             -(-1)^{\kt}\dkl \L   \label{euler}
\end{equation}
for Lagrangians in question.
\begin{proposition}
   The identity ~(\ref{euler}) is the infinitesimal form of the condition:
   \begin{equation}
      \L(p\cdot g)=\L(p)\cdot \ber g, \label{ber}
   \end{equation}
for any matrices $g$ taking values in the identity component of the
supergroup $\glpq$.
\end{proposition}
{\sc Proof. } Consider elementary
transformations  and
multiplication by $e^{\l}$ for the columns of $(\pak)$.\qed

This suggests two possible conditions for the behaviour of $\L$ under
generic transformations of the argument $p$: either they multiply it by $\ber
g$ or by $\ber_{0,1} \,g$. (Compare with two kinds of forms in the
``standard'' theory.)

\begin{definition}\label{cov}
   The Lagrangians satisfying ~(\ref{ber})
(or the similar condition with $\ber_{0,1}$) for generic $g$ are called {\it
covariant}. The cases of $\ber$ and $\ber_{0,1}$ will be distinguished
as covariance of the {\it first} or {\it second kind}.
\end{definition}

For $p|q\leqslant n|m$ the covariant dual Lagrangians
suit for integration over
$p|q$-co\-dimen\-sional submanifolds. If a cooriented submanifold
$P\subset M$ is specified by
(nondegenerate) copaths $f_{\a}$ for a covering $(U_{\a})$, then one can
define an integral over $P$ as the sum of
integrals over all $f_{\a}$, using a partition of unity.
This is well-defined due to covariance.
Actually, the covariant Lagrangians of the first kind are good for
$(+|-)$-coorientation (see~\cite{git})
and those of the second kind for $(-|-)$-coorientation.
If one prefers to consider oriented submanifolds rather than
cooriented, then it is necessary to take dual Lagrangians with local
coefficients (namely, $\sign \det TM_0$ for Lagrangians of the first kind
and $\sign\,\ber TM$ for the second kind).

The covariant Lagrangians are well-defined only, if
the rank of the matrix $(\pak)_{\kt=1}$
equals $0|q$.
In particular, the odd codegree must be less than the
odd dimension of the supermanifold:\/ $q\leqslant m$.

{}From now on we demand that the odd equations of any copath are
independent, i.e. the set of gradients $df^K$, $\kt=1$, has rank $0|q$.
At the same time, the independence of even equations
is not required. On $M^{n|m}$ their number
may very well be greater than $n$. Such copaths
do not correspond to submanifolds.
For dependent even $f^K$ the integral ~(\ref{S}) is senseless,
since the delta-function $\d(f)$ is not well-defined.

\begin{example}
   Let $\s=\s(x,\theta)Dx$ be a Bernstein--Leites integral form with local
coefficients in the sheaf $\sign\det TM_0$. Let $\s(x,\theta)$ be
a polynomial in $\theta_A$ of degree $k$.
Set
\begin{equation}
    \L(x,p)=\int\limits_{{\Bbb R}^{0|p}}\!\! D\ph\,\, \s(x,p\cdot\ph),
\end{equation}
which is well-defined for any number $p$.
The integral vanishes for $p\neq k$. For $p=k$ this integral transform
establishes an isomorphism between ``exterior polynomials'' $\s$
and skew-symmetric
multilinear functions of $p^1, \dots, p^k$. Note that if $m>0$,
these functions do not vanish for $k>n$.
\end{example}

We see that the co\-variant dual
Lagrangians of co\-degree $p|0$ are
in $1-1$-corres\-pondence with the Bernstein--Leites integral
forms of degree $n-p$. In particular, we obtain an example of
non-vanishing Lagrangians of negative even degree $n-p$.
Although the
Lagrangians of negative degree cannot be integrated
(like
integral forms of negative degree),
they constitute
an important part of the theory.

\begin{theorem}
    For covariant Lagrangians $\L$ the Euler--Lagrange equations for the
action ~(\ref{S}) are equivalent to
\begin{equation}
    (-1)^{\at\kt}\der{}{x^A}\left(\der{\L}{\pak}(x,\der{f}{x})\right)
          \equiv 0\mod (f)        \label{covel}
\end{equation}
(i.e., the left-hand side vanishes on the surface $f^K=0$).
\end{theorem}
{\sc Proof. } Expand the equation~(\ref{eulagr}) and
substitute ~(\ref{euler}) into it.
Then the terms with the derivative of the delta-function
cancel and we obtain
\begin{equation}
  \d S=\pm\int Dx\,\yk \,\d(f)(-1)^{\dt\kt}\left((-1)^{\kt\at}
                  \der{}{x^A}\left(\der{\L}{\pak}(x,\der{f}{x})\right)
      \right), \label{var}
\end{equation}
for an arbitrary finite variation $Y$. Hence if $\d S=0$, the expression
in the outer brackets annihilates $\d(f)$. That means that it belongs to
the ideal generated by $f^K$.\qed

\begin{definition}
    We call a covariant dual Lagrangian $\L$ {\it closed}, if
    \begin{equation}
    (-1)^{\at\kt}\der{}{x^A}\left(\der{\L}{\pak}(x,\der{f}{x})\right)=0
     \label{hovik}
    \end{equation}
for any $K$.
\end{definition}

That means that the action~(\ref{S}), if it makes sense,
is identically stationary for any copath.
(The equation ~(\ref{hovik}) is due to O.M.Khudaverdian~\cite{hovik}.
He introduced it
for ``dual densities'' in the sense
of~\cite{gajduk,hovik}.
)

Following the analogy with the ``standard'' theory of the previous section,
one may wish to define a differential for dual Lagrangians
(of degree one w.r.t.
Definition~\ref{deg}),
on the basis of the variation formula~(\ref{var}).
Advancing the degree of $\L$ (or the formal dimension of the copaths),
in our dual
language means reducing the number of functions $f^K$. Unfortunately, this
is not suggested by ~(\ref{var}). But, euristically, variation of any
object means introducing an additional degree of freedom, corresponding to
the variation parameter. This suggests an approach to the formula
{}~(\ref{var}) by enlarging the original ``configuration space'' $M^{n|m}$.
To do so, we are going now to
introduce a sort of  ``mixed'' theory, combining
the ``dual'' approach with the theory
of the Lagrangians of paths of the previous section.
%
%
%
%
%
%
\section {Mixed forms.}
Consider a direct product $M^{n|m}\times{\Bbb R}^{r|s}$ for some ${r|s}$,
with direct product coordinates on it.
Let $t^F$ denote coordinates on ${\Bbb R}^{r|s}$. Consider copaths
$(f^K)$ on
$\mbig$.
Two particular cases are $(x^A-x^A(t))$,
which are equivalent to paths in $M$,
and $(t-t(x),f^*(x))$, which correspond
to copaths $f^*(x)$ on $M$.
As above, in a formal way consider an integral
\begin{equation}
   S[f]=\int\limits_{\ubig}\!\!Dx\,Dt\,\d(f)\,\L\left(x,\der{f}{x},\der{f}{t}
                                            \right).
      \label{SS}
\end{equation}

There is a slight but important difference between the integrals~(\ref{SS})
and~(\ref{S})
with $M$ substituted by $M\times\Bbb R^{r|s}$.
Here in ~(\ref{SS}) the function $\,\L(x,p,w)$ does not depend on $t^F$.
That means that $\L$
cannot be treated as a geometrical object (a dual Lagrangian) w.r.t.
changes of both $x^A$, $t^F$. Instead, now the functional $S$
also depends on the coordinate
system on $\Bbb R^{r|s}$. We shall demand that $\L$ behaves as dual
Lagrangian w.r.t. $(x^A)$. So the integral ~(\ref{SS}) does not depend on
the choice of coordinates on $M$. We call functions $\L$ {\it mixed
Lagrangians\,} on $M$. The number $p|q$ is called a {\it codegree:\/}
$\codeg\L$, and $r|s$ is called an {\it additional degree:\/} $\addeg\L$.
The {\it degree\/} of $\L$ is naturally calculated as
\begin{equation}
\begin{aligned} \label{mdeg}
    \deg\L=\dim M +\addeg\L -\codeg\L\\=n+r-p\,|\,m+s-q.
\end{aligned}
\end{equation}
(As before, we consider  Lagrangians of any degree,
although the integral ~(\ref{SS}) has sense only for
$\deg\L\geqslant 0$).

\begin{example} Let a $(-|-)$-oriented closed submanifold
$P\subset M$ be specified
(locally) by
a ``mixed'' system of equations $f^K(x,t)=0$. Then for
a Bernstein--Leites pseudodifferential form $\o=\o(x,dx)$
\begin{equation}
    \int_P \o=\pm\int Dx\,D(dx)\,Dt\,D(dt)\,\d(f)\,\d(df)\, \o(x,dx).
\label{pdf}
\end{equation}
Indeed, the right-hand side is independent on the
choices of coordinates and equations (compatible with the orientation) and
coincides with the left-hand side after elimination of $t$.
\end{example}
\begin{example}
    Analoguosly, for a
(pseudo)integral Bernstein--Leites form
$\s=\s(x,\theta)Dx$  with local coefficients in the ``orientation sheaf''
$\sign\det TM_0$
\begin{equation}
    \int_P \s=\pm\int Dx\,Dt\,D\ph\,\d(f)\,
              \d\left(\der{f}{t}\ph\right)\,
                \s\left(x,\der{f}{x}\ph \right),
\label{pif}
\end{equation}
where $\ph$ runs over $\Bbb R^{q|p}$.
(A possible proof is to reduce~(\ref{pif}) to ~(\ref{pdf}) by
Fourier--Hodge integral expansion of $\s$.)
\end{example}

We stress that there is no analogue of the integrals
{}~(\ref{pdf}) and~(\ref{pif}) in the
case of $(+|-)$-orientation.

\medskip
The class of arbitrary mixed Lagrangians is too wide.
Indeed, if we specify a
submanifold by mixed copaths, the integral~(\ref{SS}) should not
depend neither on coordinates on $\Bbb R^{r|s}$
nor on the choice of equations (up to signs related to
orientation).
And if we eliminate $t^F$ from $(f^K(x,t))$, obtaining a copath of the form
$(t-t(x), f^*(x))$ on $\mbig$, then only
$f^*(x)$ is
relevant for any reasonable ``integration object that lives on $M$''.
Naturally, the mixed Lagrangians supplied by model examples above
meet these additional requirements.

\begin{definition}\label{matr}
    A mixed Lagrangian $\L=\L(x,p,w)$ is
{\it right-covariant},\, if
\vspace{-1ex}
\begin{equation}\label{rcov}
    \L(x,pg,wg)=\L(x,p,w)\cdot\ber g,
\end{equation}
is {\it left-covariant},\, if
\vspace{-1ex}
\begin{equation}\label{lcov}
    \L(x,p,hw)=\ber h \cdot \L(x,p,w)
\end{equation}
and is {\it admissible,}\, if
\vspace{-1ex}
\begin{equation} \label{adm}
    \L(x,p+aw,w)=\L(x,p,w).
\end{equation}
Here $g\in\gl(p|q)$, $h\in\gl(r|s)$, $a\in\mat(n|m\times r|s)$ are
arbitrary. (As in Definitions~\ref{cov1},~\ref{cov}, variants of
{}~(\ref{rcov}) and ~(\ref{lcov}) with $\ber$ or $\ber_{0,1}$
are possible and they will be distinguished, where necessary,
in the same way, by the words ``first'' or ``second'' kind.)
\end{definition}

The infinitesimal versions can be written for all
{}~(\ref{rcov}),~(\ref{lcov}),~(\ref{adm}). For example, the equation
\begin{equation}
\wfk\der{\L}{\pak}=0 \label{admi}
\end{equation}
is equivalent to the admis\-sibility.

\begin{proposition}
    If $\L$ is admis\-sible, then the action $S[f]$ of the copath
$(t-t(x),f^*(x))$ does not depend on the function $t(x)$.
\end{proposition}
\vspace{-1ex}
{\sc Proof. }Substitute
$\pag\!=-\derl{t^G}{x^A}(x)$,
$\pak\!=\derl{f^{*K}}{x^A}(x)$, $\wfg\!=\dfg$, $\wfk\!=0$, where the
superscript $K$
runs over the equations ``on $M$''.
Using the invariance property~(\ref{adm})%
, we can eliminate the dependence on $\derl{t^G}{x^A}(x)$
without altering $\pak$.\qed

\medskip
Now let's turn to the variation of the action ~(\ref{SS}). If $\L$
is right-covariant, then
acting formally as above for~(\ref{S}) we immediately obtain that
\begin{multline}
  \d S=\pm
   \int\! Dx\,Dt\,\d(f)\,\yk
      \left(
            (-1)^{\kt\at}\der{}{x^A}\left(\der{\L}{\pak}\right) +
            (-1)^{\kt\ft}\der{}{t^F}\left(\der{\L}{\wfk}\right)
      \right). \label{varr}
\end{multline}
Note that ~(\ref{varr}),~(\ref{var}) both include the
second derivatives of $\fk$.
The condition for
$\d S$ not to depend on second derivatives is the following equations:
\begin{align}
\secder{\L}{\pak}{\pbl}+(-1)^{\at\bt+(\at+\bt)\lt}\secder{\L}{\pbk}{\pal}&=0,
\label{eq1}\\
\secder{\L}{\pak}{\wfl}+(-1)^{\at\ft+(\at+\ft)\lt}\secder{\L}{\wfk}{\pal}&=0,
\label{eq2}\\
\secder{\L}{\wfk}{\wgl}+(-1)^{\ft\gt+(\ft+\gt)\lt}\secder{\L}{\wgk}{\wfl}&=0,
\label{eq3}
\end{align}
which are analoguous to~(\ref{bas1}). Equation~(\ref{eq1}) was suggested
in~\cite{hovik} for ``dual densities''.
Like for ~(\ref{bas1}), the odd-odd
part of these equations reminds equations that are
familiar in classical integral
geometry~\cite{john,ggg}.
Actually, from the abstract point of view all the equations~(\ref{bas1}),
(\ref{eq1}), (\ref{eq2}), (\ref{eq3}) have the same structure. They are
distinguished only by the notation for independent variables.
We shall call equations~(\ref{bas1}),~(\ref{eq1}--\ref{eq3})
and similar the {\it fundamental equations}.\,
As we shall show in a separate paper, these equations
have a homological interpretation.

Let us discuss those properties of mixed Lagrangians which are connected
with the elimination of variables from copaths. This may be treated as a
pure algebraic problem concerning abstract
functions of matrices with given conditions on
their dependence on columns and rows. Consider a matrix with a chosen row and
a chosen column of the same parity $\a$. Denote their common element by $u$,
other elements of a chosen row by $\wk$ and of a chosen column by $\va$, and
all other matrix entries by $\pak$.
The fundamental equations and various properties from Definition~\ref{matr}
have natural sense for functions of $p,w,v,u$. We'll refine only
the notion of admissibility, saying that a function of a matrix is
{\it admissible w.r.t. a given row,\/} if it is invariant under
elementary transformations of other rows by the given one.

\begin{lemma} \label{exc}
Take $\L=\L(p,w,v,u)$ and define the function
$\L^*=\L^*(p,v)$ by the formula
\begin{equation} \label{tuda}
    \L^*(p,v):=\L(p,0,v,1).
\end{equation}
Conversely, take $\L^*=\L^*(p,v)$ and define $\L=\L(p,w,v,u)$ as follows:
\begin{equation}\label{obratno}
    \L=\L(p,w,v,u):=u^{(-1)^{\a}}\L^*(p-vw/u,\,v/u).
\end{equation}
In the domain where $u$ is invertible the equations~(\ref{tuda})
and~(\ref{obratno}) supply mutually inverse isomorphisms of the spaces of
right-covariant functions of matrices $(p,w,v,u)$ and $(p,v)$.
For right-covariant $\L$ and $\L^*$ the
following properties of  $\L$ are equivalent to the properties of $\L^*$:

{\em (a) } $\L$ is admissible and homogeneous
of degree $(-1)^{\a}$
w.r.t. the chosen row $(w,u)$.  $\L^*$ does not depend on $v$.

{\em (b) }$\L$ is left-covariant for some of the rows $(p_A, \va)$.
$\L^*$ is left-covariant for the same rows.

{\em (c) }$\L$ is admissible w.r.t. a row $(p_A, \va)$.
$\L^*$ is admissible w.r.t. the same row.

{\em (d) }$\L$ satisfies the fundamental equations w.r.t. the variables
$\pak\!,\,\wk\!,\,\va,\,u$. $\L^*$ satisfies the fundamental equations
w.r.t. $\pak\!,\,\va$.
\end{lemma}

We omit the proof of Lemma~\ref{exc} because of the lack of room.
Of course, nontrivial are
the ``converse'' statements, i.e. passing back from $\L^*$ to $\L$.
Although our
proof includes rather long calculations, they seem quite beautiful.
Especially amazing is how the fundamental equations reproduce
themselves
after elimination of some independent variables using
covariance conditions. (We have discovered this fact experimentally,
in special cases, about ten
years ago, investigating the system ~(\ref{bas1}).)

{\sc Remarks. }1. The formula ~(\ref{obratno}) is written under assumption
of right-covariance of the
first kind. In the case of the second kind it is necessary to replace
$u^{(-1)^{\a}}=u^{-1}$ (for $\a=1$)
in~(\ref{obratno}) by $|u|^{-1}$.

2. Suppose $\L=\L(p,w)$ is admissible. It follows from Lemma~\ref{exc} that
right- and left-covariance conditions are
compatible for $\L$, only if they are
of the same kind. That means that ~(\ref{rcov}) and ~(\ref{lcov}) may
include either $\ber g$ and $\ber h$ or $\ber_{\!0,1}\,g$
and $\ber_{\!0,1}\,h$
but no other combination. (Else we come to a contradiction.)

\medskip
\begin{definition} \label{mixed}
    We call a right-covariant Lagrangian $\L=\L(x,p,w)$ a {\it mixed
form\/} on $M$, if it is left-covariant and admissible and
satisfies the equations~(\ref{eq1}--\ref{eq3}).
If $r|s=0$, then
$\L=\L(x,p)$ is called a {\it dual form}.
\end{definition}

Naturally distinguished are mixed forms of first and second kind.
Since for a mixed form $\L(p,w)$ both the odd momenta $(p^K,w^K)$
and the odd ``additional velocities'' $w_F$
must be linear independent, the following inequality holds:
\begin{equation}
   s\leqslant q\leqslant m+s,
\end{equation}
where $\codeg \L=p|q$, $\addeg\L=r|s$, $\dim M=n|m$. At the same time,
the even degree of $\L$, which is $n+r-p$, can be any integer, positive or
negative.

\begin{theorem}\label{iso}
    {\em (a) } Let $k,l\geqslant 0$. The space of mixed forms
of additional
degree $r+k\,|\,s+l$ and codegree $p+k\,|\,q+l$ is naturally isomorphic to
the space of  mixed forms of additional degree $r|s$ and codegree $p|q$.

{\em (b) } The space of mixed forms of additional degree $r|s$ and
codegree $n|m$, where $n|m=\dim M$, is naturally isomorhic to the space of
forms in the sense of Definition \ref{form} of degree $r|s$, ``twisted'' by
local coefficients $\sign\det TM_0$ (for the first kind) or $\sign \ber TM$
(for the second kind).
\end{theorem}
{\sc Proof. }
(a)
Consider a mixed form $\L$ of codegree $p+k\,|\,q+l$ and additional degree
$r+k\,|\,s+l$.
We shall write the arguments of forms as matrices.
Let $\left(\begin{array}{c}p \\w\end{array}\right)=
\left(\begin{array}{cc}
         \pone & \ptwo  \\
         \woo & \wot\\
         \wto & \wtt
\end{array}
\right)$. Here we distinguish the last $k|l$ of
$p+k\,|\,q+l$ momenta (written as columns) and the last $k|l$ of $r+k\,|\,s+l$
``additional velocities'' (the last $k|l$ rows).
For $\L=
\L\left(\begin{array}{c}p \\w\end{array}\right)$, set $\L^*
\left(\begin{array}{c}p^* \\w^*\end{array}\right):=\L
\left(\begin{array}{cc}
         p^* & 0  \\
         w^* & 0\\
         0 & 1
\end{array}\right)$. Obviously, $\L^*$ is a mixed form
of codegree $p|q$ and additional degree $r|s$. In the open domain where the
submatrix $\wtt$ is invertible, the form $\L$ can be expressed in terms of
$\L^*$ as follows:
\begin{equation}\label{inv}
    \L\left(\begin{array}{c}p \\w\end{array}\right)=
    \L^*\left(\begin{array}{c}
\pone-\ptwo\wtt^{-1}\wto \\
\woo-\wot\wtt^{-1}\wto \end{array}\right)\cdot\ber'\wtt.
\end{equation}
Here $\ber'$ stands for $\ber$ or $\ber_{0,1}$ for forms of first and
second kind respectively. (Use right-covariance and then admissibility or
left-covariance and then admissibility and right-covariance.) Consider the
formula ~(\ref{inv}) for arbitrary mixed forms $\L^*$.
The transformations of the arguments that we consider
can be reduced to steps, to any of which Lemma~\ref{exc} is
applicable.
Hence, the function $\L$, defined in the domain where $\wtt^{-1}$
exists, possesses all properties of a mixed form. It can be shown that $\L$
uniquely extends to all possible values of its arguments. First, by the
skew-symmetry w.r.t. both the momenta and the additional velocities of
given parity (which follows from the covariance conditions),
$\L$ uniquely extends  to the domain where {\it some\/}
$(k|l\times k|l)$-submatrix of $w$, not necessarily $\wtt$,
is invertible. That is, to the domain $\/\rank
w\ge k|l$. Since the rank of the odd rows of $w$ is $s+l\ge l$, this is a
condition only on the even-even block of $w$,
which has dimension $r+k\times p+k$.
The fundamental equations imply that $\L$ is multilinear both in the
even momenta and the even additional velocities. Thus, as a polynomial,
it is uniquely
determined by its restriction to an open domain. So we have proved that the
maps $\L\mapsto\L^*$ and $\L^*\mapsto\L$ are indeed the desired
mutually inverse isomorphisms. Their coordinate independence
is checked immediately.

(b) Let $\L=\L\left(\begin{array}{c}p \\w\end{array}\right)$ be a mixed form
of codegree $n|m$.
Define  $L(\dot x)
:=\L\left(\begin{array}{c}1 \\ \dot x\end{array}\right)$. In other
coordinate system it would be $L'(\dot x')=L'(\dot x\,\derl{x'}{x})=
\L'\left(\begin{array}{c}1 \\ \dot x\,\derl{x'}{x}\end{array}\right)=
\L\left(\begin{array}{c}\derl{x'}{x} \\ \dot x\,\derl{x'}{x}
\end{array}\right)
\cdot D_{1,0}x/D_{1,0}x'=\L\left(\begin{array}{c}1 \\
\dot x\end{array}\right)\cdot\ber'\,\derl{x'}{x}\cdot D_{1,0}x/D_{1,0}x'=\pm
L(\dot x) $, where $\pm$ equals $\sign\det(\derl{x'}{x})_{00}$ or
$\sign\,\ber\derl{x'}{x}$ for the first or the second kind, respectively.
Obviously, $L$ is covariant and satisfies the fundamental equations.
Conversely, take any form $L$, twisted by the appropriate local
coefficients, and set $\L\left(\begin{array}{c}p \\w\end{array}\right)
:=L(wp^{-1})\cdot\ber'p$, in the domain where $p$ is invertible.
{}From Lemma
{}~\ref{exc} it follows that $\L$ is a mixed form. And, as above, $\L$ uniquely
extends to all values of arguments.\qed

\medskip
The isomorphisms of Theorem~\ref{iso} are compatible with integration. For
a submanifold $P\subset M$ represented by mixed copaths, the
integral of a mixed form $\L$ over $P$, which is locally expressed by an
integral over a copath, after the elimination of any of the
parameters $t^F$ becomes an integral  of the corresponding
$\L^*$ over the new copath of less codimension.
The same is true if one eliminates the coordinates $x^A$ in order to get
an explicit local parametrization of $P\subset M$.
Then the integral of $\L$ over $P$ becomes an integral
of the corresponding $L$.

\medskip
Now we are ready to introduce a differential for mixed forms.
First, let
$\L$ be an arbitrary Lagrangian of codegree $p|q$ and additonal degree
$r|s$.
\begin{definition}
   Define its {\it differential\/} $\D\L$ as follows:
   \begin{equation}\label{dif}
      \D\L:=(-1)^{r}w_{r+1}{\!\!}^K
        \left(
            (-1)^{\at\kt}\der{}{x^A}\left(\der{\L}{\pak}\right) +
            (-1)^{\ft\kt}\der{}{t^F}\left(\der{\L}{\wfk}\right)
         \right).
   \end{equation}
(The sign $(-1)^{r}$ is inserted for convenience.)
\end{definition}
Here $\D\L$
depends on extra even ``additional velocity''  $w_{r+1}$.
Since the formula ~(\ref{dif}) includes the total differentiation
w.r.t. $x^A$ and $t^F$, the function $\D\L$ also depends on
an array of new variables $\ddot
f_{AB}$, $\ddot f_{AF}$, $\ddot f_{FG}$, unless ~(\ref{eq1}--\ref{eq3}) are
satisfied.

If $\L$ is a form, the formula for $\D\L$ is simplified:
\begin{equation} \label{diff}
      \D\L=(-1)^{r}w_{r+1}{\!\!}^K
            (-1)^{\at\kt}\der{}{x^A}\der{\L}{\pak}
\end{equation}
(because there is no explicit dependence on $t$).

\begin{theorem}\label{main}
    The definition of $\D$ does not depend on the choice of coordinates.
If $\L$ is a form, so is $\D\L$.
The identity $\D\,^2=0$ holds. The differential $\D$ commutes with the
isomophisms {\em (a)} and {\em (b)} of Theorem~\ref{iso}.
\end{theorem}
{\sc Proof of Theorem \ref{main}. }
Consider a change of coordinates on $M$: $x^A=x^A(x')$.
Since $\L(x,p,w)$ is a volume form component (a ``scalar density of tensor
weight one''), its derivative w.r.t. a component of a momentum
${\frak X}_K^A=(-1)^%
{\at(\tilde{\L}+1)}\derl{\L}%
{\pak}$ behaves as a vector density of weight one (for any $K$).
So the divergence  $\div{\frak X_K}
=(-1)^{\at(\tilde {\frak X}_K + 1)}\derl{\frak X_K^A}{x^A}$
is a well-defined scalar density. This is exactly the
first term inside the brackets in ~(\ref{dif}). As for the second term, it's
tensor behaviour is obviously the same as that of $\L$.
Thus we proved that
the operation $\D$ is well-defined for any Lagrangian $\L$.
It is convenient to divide the remaining proof into the following Lemmas.
Let $\L$ be a mixed form. Note that below the indices $F$, $G$ and
similar
enumerate  the ``old'' additional velocities $w_F$, which do not include
$w_{r+1}$.
\begin{indentedtext}
\begin{lemma}\label{R}
     $\D\L$ is right-covariant.
\end{lemma}
{\sc Proof. }Consider the linear transformations of the momenta
$(\pak,\wfk,w_{r+1}^K)$, the arguments of $\D\L$. The components
$(w_{r+1}^K)$ and the partial derivatives $\derl{}{\pak}$ transform
contragrediently. So the transformations reduce to
transformations in the argument of the form
$\L$ and we can apply the right-covariance condition.\qed
\begin{lemma}\label{L}
     $\D\L$ is left-covariant.
\end{lemma}
{\sc Proof. }Consider the linear transformations of the ``velocities''
$(w_F, w_{r+1})$. The transformations among $(w_F)$ obviously commute with
$\secderl{}{x^A}{\pak}$, so we can apply the left-covariance of $\L$.
The homogeneity  in $w_{r+1}$ is evident. Consider an
elementary transformation of the form $w_{r+1}\mapsto w_{r+1}+\e\, w_F$.
It leaves $\D\L$ invariant due to the admissibility of $\L$.
Finally, consider the elementary transformation $w_F\mapsto w_F+\e\, w_{r+1}$.
It leaves $\D\L$ invariant because of the equation~(\ref{eq2}).\qed
\begin{lemma}\label{A}
     $\D\L$ is admissible.
\end{lemma}
{\sc Proof. }Consider an elementary transformation of the row $p_A$
by the row $w_F$. It obviously commutes with $\D$ and we can apply the
admissibility of $\L$. Now consider elementary transformations
by $w_{r+1}$.
They leave $\D\L$ invariant because of the equation~(\ref{eq1}).\qed
\end{indentedtext}
\medskip
Since for a mixed form $\L$ the Lagrangian $\D\L$ does not depend on second
derivatives, we may again apply to it the operation $\D$
by the formula~(\ref{dif}) and obtain
a mixed Lagrangian $\D\D\L$.
\begin{indentedtext}
\begin{lemma}\label{ddva}
     $\D\D\L=0$.
\end{lemma}
{\sc Proof. }We need to calculate the expression
\begin{multline} \label{dd}
\qquad w_{r+2}^L\left((-1)^{\bt\lt}\der{}{x^B}\der{}{\pbl}
+(-1)^{\gt\lt}\der{}{t^G}\der{}{\wgl}\right. \\
\left.+\der{}{t^{r+1}}\der{}{w_{r+1}^K}\right)
\left(w_{r+1}^{K}(-1)^{\at\kt}\der{}{x^A}\der{\L}{\pak}\right),
\end{multline}
where the differentiation in the outer brackets is ``total''. By
straightforward but rather tedious calculation, (\ref{dd}) reduces to
a sum of six terms, each of which vanishes by
the equations~(\ref{eq1}--\ref{eq3}).\qed
\end{indentedtext}
\medskip
\noindent
Lemma~\ref{ddva} in particular
implies that $\D\L$ satisfies the fundamental equations.
Thus, by Lemmas~\ref{R},~\ref{L},~\ref{A} and~\ref{ddva}, $\D$ indeed maps
mixed forms to mixed forms and, by Lemma~\ref{ddva}, \quad $\D\,^2=0$\quad
on mixed forms.

Now let's prove that $\D$ commutes with the isomorphisms of
Theorem~\ref{iso}.
Let $\L^*$ be a mixed form of codegree $p|q$ and additional degree $r|s$.
Consider $\L(x,p,w)=\L^*(x,\pone-\ptwo\wtt^{-1}\wto, \woo-\wot\wtt^{-1}\wto)
\cdot\ber'\wtt$ (notation of Theorem~\ref{iso} except that,
for typographical
reasons, we are not consistent in writing matrices). The codegree and
additional degree for $\L$ is advanced by $k|l$. Apply the operator $\D$ to
$\L$. By immediate check, $(-1)^{r+k}\D\L$ is expressed in terms of
$(-1)^r\D\L^*$ by the same formula as $\L$ is expressed in terms of $\L^*$.
This coincides with the desired isomorphism up to the sign $(-1)^k$,
because one needs to transpose $w_{r+k+1}$ with the last $k$ even rows $w_F$
(in the argument of $\D\L$). Thus the signs cancel. The proof for the
isomorphism (b) is similar.
\smallskip
{\noindent\sc End of the Proof of Theorem~\ref{main}. }

\bigskip
Denote the space of mixed forms of first or second kind
twisted by $\sign\det TM_0 $ or $\sign\ber TM$, respectively,
by $\O^{r|s}_{p|q}$.
Since the author has
problems with drawing in $4$-dimensional space, below the numbers $q$
and $s$ are fixed.
{}From Theorem~\ref{main} we obtain the following commutative diagram.

$$
   \renewcommand{\om}[2]{{\O^{#1|s}_{#2|q}}}
   \renewcommand{\to}{@>{\D}>>}
   \newcommand{\up}{@AAA}
   \begin{CD}
    {\ldots}  @. {\ldots} @. {\ldots} @. {\ldots} @. {\ldots} @.{\ldots}\\
     @.             \up          \up         \up         \up         \up\\
     0     @>>>   \om{0}{2} \to\om{1}{2}\to\om{2}{2}\to\om{3}{2}\to\om{4}{2}
                                                                \to{\ldots}\\
     @.              @.           \up        \up          \up      \up \\
     { }      @.     0    @>>>  \om{0}{1}\to\om{1}{1}\to\om{2}{1}\to\om{3}{1}
                                                                \to{\ldots}\\
     @.              @.           @.          \up         \up       \up \\
     {}       @.     {}   @.      0 @>>>    \om{0}{0}\to\om{1}{0}\to\om{2}{1}
                                                                \to{\ldots}\\
     {\ldots} @.{n-2}        @.  {n-1} @.     {n}     @. {n+1}  @.
                                                        {n+2}
                                                               @. {\ldots}
   \end{CD}
$$
Here the vertical arrows are the isomorphisms from Theorem~\ref{iso}.
The numbers in the line below the diagram are equal to the even
degrees of mixed
forms in the corresponding column. (The odd degree $m+s-q$ is the same
for the whole diagram.)
Let's denote mixed forms without any twist by $\si{r}{s}{p}{q}$. The same
diagram can be drawn for them, too.
It is natural to pass to a direct limit (which is very simple here).

\begin{definition} Let $r\in\Bbb Z$ be any integer and let
$0\le s\le m$. Define the spaces of ``stable'' forms by
\newcommand{\dlim}{{\mathop{\varinjlim}\limits_{N,M}}}
   \begin{equation}
     \begin{aligned}
       \Om{r}{s}:=\dlim\om{N}{M}{n-r+N}{m-s+M},\\
       \Si{r}{s}:=\dlim\si{N}{M}{n-r+N}{m-s+M}.
     \end{aligned}
   \end{equation}
\end{definition}
Obviously,
\begin{equation}
     \begin{aligned}
       \Om{r}{s}=\om{N}{M}{n-r+N}{m-s+M},\\
       \Si{r}{s}=\si{N}{M}{n-r+N}{m-s+M},
     \end{aligned}
\end{equation}
for $N\ge 0$, $N\ge r-n$, $M\ge 0$, $M\ge s-m$. In particular,
if $r\ge 0$, then
\begin{equation}
    \Om{r}{s}=\om{r}{s}{n}{m}=\O^{r|s}
\end{equation}
(forms in the sense of \S~1). And if $r\le n$, then
\begin{equation}
    \Om{r}{s}=\om{0}{0}{n-r}{m-s}
\end{equation}
(twisted dual forms). Thus the ``stable complexes'' $\Om{\cdot}{s}$ can be
completely described in terms of $\O^{\cdot|s}$  and
$\O_{\cdot|q}:=\om{0}{0}{\cdot}{q}$. The differential on $\Om{\cdot}{s}$
induce the following differential for dual forms:
\begin{equation}
    \DD\L:=(-1)^{k-1}\der{}{x^A}\der{\L}{p_A{}^{k}},
\end{equation}
for $\L\in\O_{k|q}$. (It follows from Theorems~\ref{iso},~\ref{main}
that $\L$ is
closed, iff $\DD\L=0$.) The complexes $\O^{\cdot|s}$ and $\O_{\cdot|m-s}$
are sewed together as follows:

\advance\leftskip by-1.5cm
{\parbox{0em}%
{$$
    \renewcommand{\a}[1]{{\O^{#1|s}%
                              }}
    \renewcommand{\b}[1]{{\O_{#1|m-s}%
                              }}
    \renewcommand{\to}{@>>>}
    \newcommand{\tot}{@>{\D}>>}
    \newcommand{\tob}{@>{\DD}>>}
    \begin{CD}{}@.{0}\mqquad\to{\!\!\mquad\a{0}}\mquad\,\tot\mqquad
     \a{1}\mquad\,\tot{\ldots}\tot\!\!\mquad\a{n}\mquad\,\tot\a{n+1}\tot
                                                                 {\ldots}\\
     @.  @. \!\!\mquad @AAA\mqquad @AAA  @.\!\!\mquad @AAA  @.\\
          {\ldots}\tob\b{n+1}\mquad\,{\larrow{-7}{10}}\tob\b{n}\mquad\,\tob
     \b{n-1}\mquad\,\tob{\ldots}\tob\b{0}\mquad\,\rarrow{-5}{10}\to\mqquad{0}
    \end{CD}
$$}

\advance\leftskip by1.5cm
Here the vertical arrows are isomorphisms and the dashed arrows are defined
so to make the diagram commutative.
Recall from \S 1 that for forms as Lagrangians on paths
it was necessary to introduce an augmentation. It can be shown that
the ``exact $0|s$-forms'' $B^{0|s}$ introduced in a formal fashion in \S 1
coincide with the image of $\h:\O_{n+1|m-s}\to\O^{0|s}$.

%
%
%
%
%
%

%
%
\end{document}